\documentclass[12pt,A4]{article}
\usepackage{amssymb, amsmath, enumerate, setspace}

\begin{document}
	\newtheorem{thm}{Theorem}
	\newtheorem{lm}{Lemma}
	\newtheorem{prop}{Proposition}
	\newtheorem{cor}{Corollary}
	\newtheorem{defi}{Definition}
	\newtheorem{eg}{Example}
	\newtheorem{coj}{Conjecture}
	\newtheorem{rem}{Remark}
	
	\numberwithin{equation}{section}
	
	\newcommand{\B}{{\bf Proof\ }}
	\newcommand{\A}{{$\blacksquare$\bigskip }}

	\newcommand{\CC}{\mathcal C}
	\newcommand{\CF}{\mathcal F}
	\newcommand{\CR}{\mathcal R}
	\newcommand{\cs}{\mathcal S}
	\newcommand{\CV}{\mathcal V}
	\newcommand{\ch}{\mathcal H}
	\newcommand{\ca}{\mathcal A}
	\newcommand{\cb}{\mathcal B}
	\newcommand{\cq}{\mathcal Q}
	\newcommand{\X}{\mathcal X}
	\newcommand{\W}{\mathcal W}
	\newcommand{\cu}{\mathcal U}
	\newcommand{\ic}{\mathcal I}
	\newcommand{\M}{\mathcal M}
	\newcommand{\mco}{\mathcal O}
	\newcommand{\cl}{\mathcal L}
	\newcommand{\m}{\underline{M} }
	\newcommand{\mb}{\bar{M} }
	\newcommand{\C}{\Bbb C}
	\newcommand{\N}{\Bbb N}
	\newcommand{\bbm}{\mathbb{M}}
	\newcommand{\pr}{\mathbb{P}}
	\newcommand{\R}{\mathbb{R}}
	\newcommand{\Q}{\mathbb{Q}}
	\newcommand{\s}{\mathbb{S}}
	\newcommand{\Z}{\mathbb{Z}}
	\newcommand{\I}{\mathbb{I}}
	\newcommand{\T}{\mathbb{T}}
	\newcommand{\D}{\mathbb{D}}
	\newcommand{\U}{\mathbb{U}}
	
	\newcommand{\dia}{\mathrm{dia}}
	\newcommand{\id}{\mathrm{id}}
	\newcommand{\1}{\mathbf{1}}
	\newcommand{\mo}{{}^{\scriptscriptstyle -1}}
	\renewcommand{\ij}{_{ij}}
	\newcommand{\kj}{_{kj}}
	
	\newcommand{\rf}{{}_{\scriptscriptstyle(F)}}
	\newcommand{\lf}{_{}{\scriptscriptstyle F}}
	\newcommand{\eps}{\epsilon}
	\newcommand{\Del}{\Delta}
	\newcommand{\del}{\delta}
	\newcommand{\Lan}{\Lambda}
	\newcommand{\Gam}{\Gamma}
	\newcommand{\gam}{\gamma}
	\newcommand{\ten}{\bigotimes}
	\newcommand{\op}{\oplus}
	\newcommand{\im}{\imath}
	\newcommand{\st}{\otimes}
	\newcommand{\fl}{\rightarrow}
	\newcommand{\mpt}{\mapsto}
	\newcommand{\ti}{\times}
	
	\newcommand{\Om}{\Omega}
	\newcommand{\sfi}{\varphi}
	\newcommand{\om}{\omega}
	\newcommand{\lan}{\lambda}
	\newcommand{\sig}{\sigma}
	
	\newcommand{\al}{\alpha}
	\newcommand{\be}{\beta}
	\newcommand{\te}{\theta}
	\newcommand{\all}{\forall}
	\newcommand{\ci}{\circ}
	\newcommand{\sm}{\shortmid}
	\newcommand{\ld}{\ldots}
	\newcommand{\cd}{\cdots}
	\newcommand{\p}{\partial}
	\newcommand{\n}{\nabla}
	\newcommand{\mt}{\mapsto}
	\newcommand{\se}{\subseteq}
	\newcommand{\fa}{\forall}
	\newcommand{\iy}{\infty}
	
	\newcommand{\eqn}[2]{\begin{equation}#2\label{#1}\end{equation}}

	\textheight 23.6cm \textwidth 16cm 
	\topmargin -.2in \headheight 0in
	\headsep 0in \oddsidemargin 0in \evensidemargin 0in \topskip 28pt
	
	\title{
		\textbf{Nonpure (Non)Commutative Analysis, Geometry and Mechanics}\\{ \small  part 1: differential and integral calculus}}
	
	
	\author{{\small Seyed Ebrahim Akrami}\\{\small Department of Mathematics, Semnan University, Semnan, Iran}\\{\small akramisa@semnan.ac.ir, https://orcid.org/
			0000-0002-8913-291X}}
	\maketitle
	\begin{abstract}
		We construct and study differential and integral calculus  on the space of states of a C*-algebra by equipping it with a formal smooth  structure. To achieve this goal we first concentrate on the space of nonpure states of a commutative C*-algebra as a guideline for the noncommutative case.
		In particular, we prove Stokes' theorem over the both commutative and noncommutative smooth Wasserstein space.
	\end{abstract}

\section{Introduction}
In 	these series of  papers we study analysis, topology, geometry and dynamics on the space of states of a $C^*$-algebra. Our plan is to put a formal smooth structure on this space. To prepare for performing this task we first concentrate on the commutative case. Namely, we first extend analysis, topology, geometry and dynamics from the space of pure states of a commutative $C^*$-algebra to the space of nonpure states of this algebra. Next, we extend the theory obtained in the commutative case to the noncommutative case.  
	
There are some studies on the space of nonpure states of a $C^*$-algebra. Connes was studied the metric aspect of the space of states of an arbitrary   spectral triple by putting a metric on this space which he showed that it  coincides with the Riemannian geodesic distance in the case of a commutative spectral triple. But unfortunately this metric space was not studied much until recently. Martinetti computed this metric for some classes of noncommutative examples \cite{M1}. Lorentzian distance formula in NCG has been studied by Franco \cite{F}. In the commutative case the space of nonpure states was very well-known from the time of Monge in the theory of \emph{optimal transport} and is called Wasserstein space and is the space of probability measures on a metric space \cite{V}. Rieffel showed that the Connes' metric coincides with the Wasserstein's metric in the commutative case \cite{R1,R2}. On the other hand  Ambrosio and Gigli \cite{AGS} studied analysis and geometry on a metric-measure space. They defined the concept of a smooth path over such space and in particular they studied analysis and geometry of the  Wasserstein space. Villani,  Lott, and Sturm studied the concept of Ricci curvature for a metric measure space in general and in particular for the space of nonpure states of a commutative  algebra, \cite{LV, V}.  Most  relevant works to our project are the works \cite{AG, GNT, L1, L2, LV} where they started to study  Riemannian geometry and dynamics on the space of nonpure states of a commutative $C^*$-algebra, i.e. on the space of Wasserstein space by putting the structure of a formal manifold or weak analogue of a differential and Riemannian structure on it, first introduced by Otto \cite{O}. The author was not aware of these results and the concepts and results obtained by him are independent of the above works and are consequence of a study on foundation of quantum mechanics. 

\textbf{Acknowledgment} Thanks to professor M. Golshani, for valuable discussions during my visit of department of physics at IPM Tehran Iran 2010 who made me familiar with Bohmian mechanics which is the origin of the present work.

\section{Nonpure Commutative Analysis}
\subsection{Nonpure Commutative Differential Calculus}\label{npcdc}
In this section, $M$ is a smooth oriented Riemannian manifold and the integration of a function $f$ over $M$ is computed with respect to the Riemannian volume form and is dented by $\int f$.
\begin{defi} 
i) Let $\s(M)\subset C_0(M)^*$ be the space of all states (positive linear functional of unit norm) of the $C^*$-algebra $C_0(M)$. We equip the space $\s(M)$ with the induced subspace weak*-topology from $C_0(M)^*$.\\ 
		ii) It is well-known  that $M$ can be embedded homeomorphically  in $\s(M)$ by evaluation, i.e. by assigning to each point $x\in M$ the pure state $\hat{x}\in\s(M), \hat{x}(f):=f(x),\fa f\in C_0(M).$\\
		iii) We define $\I(M)\subset \s(M)$ to be the subspace of all nonpure states of the following form
		\eqn{nonpurestate}{\int_M\rho(x)\hat{x},} where $\rho:M\fl[0,\iy)$ is smooth  with compact support and $\int\rho=1.$ Namely	
		\eqn{}{(\int\rho(x)\hat{x})(f)=\int\rho(x)f(x),~~~~~\fa f\in C_0(M).}
		We equip $\s(M)=M\cup\I(M)$  with the weak*-topology and call it \textbf{smooth convex closure} of $M$. 
	\end{defi}
	\begin{prop}
		i) A sequence $\rho_n\in\I(M)$ converges to $\rho\in\I(M)$ iff $\int\rho_nf\fl\int\rho f,\fa f\in C_0(M)$.\\
		ii) A sequence $\rho_n\in\I(M)$ converges to $p\in M$ iff $\int\rho_nf\fl f(p),\fa f\in C_0(M)$.\\
		iii) A sequence $p_n\in M$ converges to $p\in M$ iff $f(p_n)\fl f(p),\fa f\in C_0(M)$.\\
		iv) Never a sequence in $M$ converges to a point in $\I(M)$. Namely  the induced subspace topology from $\s(M)$ on $M$ is the own topology of $M$ and $M$ is a closed subspace.\\
		v) $M$ is the boundary of $\I(M)$ under the induced subspace topology. Namely for any given $p\in M$ there exists a sequence $\rho_n\in\I(M)$ converges to $p\in M$.\\
		vi) The space $(\s(M),$ weak*-topology$)$ is metrizable
		\eqn{}{d(\phi,\psi):=\sup\{\phi(f)-\psi(f)~|~|f(x)-f(y)|\le d(x,y), \fa x,y\in M \}} where $d$ is the Riemannian distance induced from the Riemannian  metric on $M$. Thus  the space $\s(M)$ is Hausdorff and a sequence $\rho_n\in\I(M)$ of nonpure points cannot converge to both a pure point in $M$ and to a nonpure point in $\I(M)$.
	\end{prop}

	Next, we convert the space of nonpure points $\I(M)$ to a formal manifold. Namely, we do not put a smooth atlas on this space  but we will say how to do differential and integral calculus on $\I(M)$. Alternatively, one can consider it as a metric space and then by the theory of metric analysis and geometry, \cite{AGS}, do differential calculus on this space. We briefly bring here this theory. First, let $(S,d)$ be a complete metric
	space. For any curve $\rho:(a, b)\se\R\fl S$ the limit
	\eqn{}{|\rho'|(t) := \lim_{h\fl0} \frac{d(\rho(t + h),
			\rho(t))}{|h|}} if exists, is called the \textbf{metric
		derivative} of $\rho$.
	
	For any absolutely continuous curve $\rho:(a, b)\se\R\fl S$ the
	metric derivative exists for Lebesgue-a.e. $t\in(a, b)$ and
	$d(\rho(s), \rho(t))\le\int_s^t |\rho'|(r) dr$ for any interval
	$(s, t)\se(a, b).$ In the case where $S$ is a Banach space with
	norm $\| \|$, a curve $\rho:(a,b)\fl S$ is absolutely continuous
	if and only it is differentiable in the ordinary sense for
	Lebesgue-a.e. $t\in(a, b)$ and we have $\|\rho'(t)\|=|\rho'|(t)$
	for Lebesgue-a.e. $t\in(a, b)$. Next in \cite{AGS}, this fact
	has been applied to the metric space of probability measures on a
	Hilbert space $X$. This space is equipped with the Wasserstein's metric. It is
	shown in \cite{AGS} that the class of absolutely curves
	$\rho_t$ in this metric space coincides with solutions of the
	continuity equation of physicists. More precisely, given an absolutely
	continuous curve $\rho_t$, one can find a Borel time-dependent
	velocity field $V_t : X\fl X$ such that $\|V_t\|_{L^p(\rho_t)}\le
	|\rho'|(t)$ for a.e. $t$ and the continuity equation
	$\frac{\p\rho}{\p t}+\n.(\rho V)=0$ holds. Conversely, if $\rho_t$
	solves the continuity equation for some Borel velocity field $V_t$
	with $\int_a^b \|V_t\|_{L^p(\rho_t)}dt < \iy$, then $\rho_t$ is an
	absolutely continuous curve and $\|V_t\|_{L^p(\rho_t)}\ge
	|\rho'|(t)$ for a.e. $t \in(a, b).$ As a consequence, we see that
	among all velocity fields $V_t$ which produce the same flow
	$\rho_t$, there is a unique optimal one with smallest
	$L^p(\rho_t,X)$-norm, equal to the metric derivative of $\rho_t$.
	One can  view this optimal field as the ``tangent" vector field to
	the curve $\rho_t$. In this paper we will not use the above mentioned  theory and as we said before we will put informal smooth structure on this space, but our definitions  and results are in agreement partially with the results of this theory. 
	
	We start by determining when a curve $\rho_t\in \I(M)$ is smooth and what  it's derivative is. To answer to this question we again imitate the case of space of pure states $M$ itself. We know that a curve $x(t)\in M$ is smooth if and only if there exists a vector $v(t)\in T_{x(t)}M$ such that
	\eqn{}{\frac{d}{dt}f(x(t))=df_{x(t)}(v(t))} for all smooth functions $f$ on $M$ and in fact $v(t)=dx(t)/dt$. By imitating the above fact we suggest the following definition of differentiability of a nonpure curve.
	\begin{defi}
		A nonpure curve $\rho_t\in \I(M)$ is called smooth if there exists a smooth time-dependent  vector field $V_t\in\X(M)$ such that
		\eqn{fracddtintrhof}{\frac{d}{dt}\int\rho f=\int\rho df(V)} for all smooth functions $f$ on $M$.
	\end{defi}
	\begin{thm}
		A nonpure curve $\rho_t$ is smooth if and only if there exists a smooth time-dependent  vector field $V_t\in\X(M)$ such that
		\eqn{f}{\frac{\p\rho}{\p t}+\n.(\rho V)=0.}
	\end{thm}
	\B  We have $\frac{d}{dt}\int\rho f=\int\frac{\p\rho}{\p t}f$. On the other hand 
	\begin{eqnarray}\int\rho df(V)&=&\int\n.(\rho f V)-\int f\n.(\rho V)\nonumber\\
		&=&-\int f\n.(\rho V).\nonumber\end{eqnarray}Thus (\ref{fracddtintrhof}) holds if and only if $\int\Big(\frac{\p\rho}{\p t}+\n.(\rho V)\Big)f=0, \fa f$. The later is equivalent with (\ref{f}). \A
	
	Note that the vector field $V$ is not unique. Thus   we redefine our concept of differentiability as follows.
	\begin{defi}\label{smoothnonpurecurve}
		A \textbf{smooth nonpure curve} is a pair $(\rho_t,V_t)$ where $\rho_t\in \I(M),V_t\in\X(M)$ and  (\ref{f}) holds.
	\end{defi}
	Next we define a several variable differentiable nonpure function.
	\begin{defi}\label{multi smooth np }
		A \textbf{smooth nonpure   function} from open subset
		$U\se\R^m$ into  $M$ is a system
		$(\rho,V_1,\cd,V_m)$ including a nonpure function $\rho:U\fl\I(M)$ such that   the induced function
		\eqn{}{\rho:U\ti M\fl[0,\iy)}is smooth,  together
		with smooth vector fields \eqn{}{V_j:U\ti M\fl TM} $1\le
		j\le m,$ such that \eqn{}{\frac{\p \rho}{\p u_j}+\n.(\rho V_j)=0.}
	\end{defi}
	
	A nonpure function $\rho:U\fl \I(M)$ is replacement for a pure
	function $f:U\fl M$ and the vector fields $V_i$ are replacement
	for the partial derivatives $V_i=\frac{\p f}{\p u_i}$. We have the equality
	$\frac{\p^2 f}{\p u_i\p u_j}=\frac{\p^2 f}{\p u_j\p u_i}$. Namely  \eqn{pAi}{\frac{\p V_i}{\p
			u_j}=\frac{\p V_j}{\p u_i}.}
	For nonpure functions we have the
	following result.
	\begin{thm}If $(\rho,V_1,\cd,V_m)$ is  a smooth
		nonpure  function from open region $U\se\R^m$ to $M$ then
		\eqn{}{\n.\Big(\rho(\frac{\p V_i}{\p u_j}-\frac{\p V_j}{\p
				u_i}-[V_i,V_j])\Big)=0.}Here $[V_i,V_j]$ means the Lie bracket of
		vector fields $V_i$ and $V_j$ with respect to $x$ variable.
	\end{thm}
	\B First proof. First we state an important identity which we discovered newly and we have not seen before elsewhere. For any two vector fields $V$ and $W$ and any function $f$ on $M$ we have
	\eqn{}{\n.\Big(\n.(f W)V\Big)-\n.\Big(\n.(f V)W\Big)=\n.(f[W,V]).} The proof is easy. Now using this identity  we have	
	\begin{eqnarray}
		\frac{\p^2\rho}{\p u_j\p u_i}&=&-\frac{\p}{\p u_j}\n.(\rho V_i)\nonumber\\&=&-\n.(\frac{\p\rho}{\p u_j}V_i)-\n.(\rho\frac{\p V_i}{\p u_j})\nonumber\\&=&\n.\Big(\n.(\rho V_j)V_i\Big)-\n.(\rho\frac{\p V_i}{\p u_j}).\nonumber
	\end{eqnarray}
	Thus 
	\begin{eqnarray}
		0&=&\frac{\p^2\rho}{\p u_j\p u_i}-\frac{\p^2\rho}{\p u_i\p u_j}\nonumber\\&=&\n.\Big(\n.(\rho V_j)V_i\Big)-\n.\Big(\n.(\rho V_i)V_j\Big)-\n.(\rho\frac{\p V_i}{\p u_j})+\n.(\rho\frac{\p V_j}{\p u_i})\nonumber\\&=&-\n.(\rho[V_i,V_j])-\n.(\rho\frac{\p V_i}{\p u_j})+\n.(\rho\frac{\p V_j}{\p u_i})\nonumber\\&=&\n.\Big(\rho(\frac{\p V_j}{\p u_i}-\frac{\p V_i}{\p u_j}-[V_i,V_j])\Big).\nonumber
	\end{eqnarray}
	Second proof. For any function $g$ on $M$, we have
	\begin{eqnarray}\frac{\p}{\p u_i}\int \rho g&=&\int \frac{\p \rho }{\p u_i}g\nonumber\\
		&=&-\int g\n.(\rho V_i)\nonumber\\
		&=&-\int\n(\rho gV_i)+\int \rho V_ig\nonumber\\
		&=&\int \rho V_ig\nonumber\end{eqnarray} by the divergence theorem and that $\rho$ vanishes at infinity.  Thus
	\begin{eqnarray}\frac{\p^2}{\p u_j\p u_i}\int \rho g&=&\int\frac{\p \rho }{\p u_j}V_ig+\int \rho \frac{\p V_i}{\p u_j}g\nonumber\\
		&=&-\int \n.(\rho V_j)V_ig+\int \rho \frac{\p V_i}{\p u_j}g\nonumber\\
		&=&-\int \n.(\rho (V_ig)V_j)+\int \rho V_jV_ig+\int \rho \frac{\p V_i}{\p u_j}g\nonumber\\
		&=&\int \rho (\frac{\p V_i}{\p u_j}+V_jV_i)g\nonumber\end{eqnarray}by the divergence theorem and that $\rho$ vanishes at infinity. 
	Hence 
	\begin{eqnarray}0&=&\frac{\p^2}{\p u_j\p u_i}\int \rho g-\frac{\p^2}{\p u_i\p u_j}\int \rho g\nonumber\\&=&\int \rho (\frac{\p V_i}{\p u_j}+V_jV_i-\frac{\p V_j}{\p u_i}-V_iV_j)g\nonumber\\&=&\int \rho (\frac{\p V_i}{\p u_j}-\frac{\p V_j}{\p u_i}+[V_i,V_j])g.\nonumber \end{eqnarray} Since $g$ is arbitrary, the desired identity is obtained.\A

	But in practice we need that $\frac{\p V_i}{\p u_j}-\frac{\p
		V_j}{\p u_i}-[V_i,V_j]$ vanishes. Thus we complete the definition of  several variable  differentiable nonpure function as follows.
	as following.
	\begin{defi}\label{multi smooth np }
		A \textbf{smooth nonpure   function} from open subset
		$U\se\R^m$ into  $M$ is a system
		$(\rho,V_1,\cd,V_m)$ including a nonpure function $\rho:U\fl\I(M)$ such that   the induced function
		\eqn{}{\rho:U\ti M\fl[0,\iy)}is smooth,  together
		with smooth vector fields \eqn{}{V_j:U\ti M\fl TM} $1\le
		j\le m,$ such that \eqn{}{\frac{\p \rho}{\p u_j}+\n.(\rho V_j)=0} and
		\eqn{pvipuj}{\frac{\p V_i}{\p u_j}-\frac{\p
				V_j}{\p u_i}=[V_i,V_j].}
	\end{defi}
	
	\begin{eg}
		Let $A\in M_{n\ti m}(\R)$ be a matrix of size $n\ti m$ with real entries whose columns are  $A_i,1\le i\le m,$ and let $\sig\in \I(\R^n)$. Then the pair $(\rho,V_1,\cd,V_m)$ given below is a smooth nonpure function from $\R^m$ into $\R^n$.
		\eqn{}{\rho(x,y):=\sig(y-Ax),~~~~~~V_i(x,y):=A_i.}
	\end{eg}
	\B We have $\frac{\p\rho}{\p x_i}=-\n \sig(y-Ax).A_i$ and $\n^y.(\rho(x,y)V_i(x,y))=\n^y.(\sig(y-Ax)A_i)=\n\sig(y-Ax).A_i$. Clearly the identities  (\ref{pvipuj}) hold.
	\A
	
	\begin{defi}\label{smoouncerfun}
		A \textbf{differentiable nonpure  function} from a manifold $P$
		into oriented Riemannian manifold $M$ is a pair $(\rho,V)$
		including a nonpure  function \eqn{}{\rho:P\ti M\fl[0,\iy)}which is
		differentiable too, together with a differentiable map
		\eqn{}{V:TP\ti M\fl TM,~~~~~~V(w,q)\in T_qM} $\fa w\in TP,q\in M$,
		such that $V(aw_1+w_2,q)=aV(w_1,q)+V(w_2,q),\fa a\in\R,w_i\in
		T_pP,p\in P,q\in M$ and for any coordinate system
		$u=(u_1,\cd,u_m)\in U\se\R^m$ for $P$ \eqn{pphi}{\frac{\p \rho}{\p
				u_j}+\n.(\rho V_j)=0} where $V_j:U\ti M\fl
		TM,V_j(u,q):=V(\frac{\p}{\p u_j}(u),q)\in T_qM.$ Here the
		divergence operator $\n.$ is with respect to variable $q\in M$. Moreover we assume that
		\eqn{fracpvi}{\frac{\p V_i}{\p u_j}-\frac{\p
				V_j}{\p u_i}=[V_i,V_j].}
		The map $V$ is called a \textbf{derivative} for the  nonpure
		function $\rho$.
	\end{defi}
	\begin{prop}\label{dprho}
		The  equations (\ref{pphi}) and (\ref{fracpvi}) do not depend to coordinate systems of $P$.
	\end{prop}
	\B Let $B_i(v):=V_{p}(\frac{\p}{\p v_i}(p))$ where $v$ is the
	coordinates of point $p$  then since $\frac{\p}{\p
		v_i}=\sum_j\frac{\p u_j}{\p v_i}\frac{\p}{\p u_j}$ we get
	$B_i(v)=\sum_j\frac{\p u_j}{\p v_i}V_j(u)$. Thus
	\begin{eqnarray}& &\frac{\p\rho(v,q)}{\p v_i}+\n.(\rho(v,q)B_i(v,q))
		\nonumber\\&=& \sum_j\frac{\p u_j}{\p v_i}\frac{\p\rho(u,q)}{\p
			u_j}+\n.(\rho(u,q)\sum_j\frac{\p u_j}{\p
			v_i}V_j(u,q))\nonumber\\&=&\sum_j\frac{\p u_j}{\p
			v_i}[\frac{\p\rho(u,q)}{\p u_j}+\n.(\rho(u,q)V_j(u,q))].\nonumber
	\end{eqnarray}Thus since the matrix $(\frac{\p u_j}{\p
		v_i})$ is invertible we conclude that $\frac{\p\rho(v,q)}{\p
		v_i}+\n.(\rho(v,q)B_i(v,q))=0$ if and only if
	$\frac{\p\rho(u,q)}{\p u_i}+\n.(\rho(u,q)V_i(u,q))=0.$
	
	Next,
	$$\frac{\p B_i}{\p v_j}=\sum_k\frac{\p^2u_k}{\p v_j\p
		v_i}V_k+\sum_{kl}\frac{\p u_k}{\p v_i}\frac{\p u_l}{\p
		v_j}\frac{\p V_k}{\p u_l}$$ Thus $$\frac{\p B_i}{\p v_j}-\frac{\p
		B_j}{\p v_i}=\sum_{kl}\frac{\p u_k}{\p v_i}\frac{\p u_l}{\p
		v_j}(\frac{\p V_k}{\p u_l}-\frac{\p V_l}{\p u_k}).$$ And
	$$[B_i,B_j]=\sum_{kl}\frac{\p u_k}{\p v_i}\frac{\p
		u_l}{\p v_j}[V_k,V_l].$$ Hence
	$$\frac{\p B_i}{\p u_j}-\frac{\p B_j}{\p u_i}-[B_i,B_j]=
	\sum_{kl}\frac{\p u_k}{\p v_i}\frac{\p u_l}{\p v_j}(\frac{\p
		V_k}{\p u_l}-\frac{\p V_l}{\p u_k}-[V_k,V_l]).$$ Thus since the
	matrix $(\p
	u_j/\p v_i)$ is invertible we conclude that\\
	$$\frac{\p B_i}{\p u_j}-\frac{\p B_j}{\p
		u_i}-[B_i,B_j]=0$$ if and only if $$\frac{\p V_k}{\p u_l}-\frac{\p
		V_l}{\p u_k}-[V_k,V_l]=0.$$ \A
\subsection{Nonpure Commutative Integral Calculus}
	In the pure commutative case when we have a smooth map $F:Q\fl M$ from a $k$-dimensional oriented manifold $Q$ to $M$ then we can pull back any $k$-differential form $\om$ on $M$ to a $k$-dimensional differential form $F^*\om$ on $Q$ and then integrate over $Q$ as $\int_QF^*\om$. Then we have the Stokes' theorem  on map $F$
	\eqn{}{\int_QF^*d\om=\int_{\p Q}F^*\om.}Proof is easy: $\int_QF^*d\om=\int_QdF^*\om=\int_{\p Q}F^*\om$ by the Stokes' theorem on $Q$. This theorem is used in the construction of pairing between singular homology and de Rham cohomology of $M$ in the de Rham theorem. Now we are going to proof a nonpure version of this theorem. 
	\begin{defi}\label{npcpull}
		Let $F=(\rho,V_1,\cd,V_k)$ be a nonpure map from $Q$ to $M$. The pullback of $k$-differential form $\om$ on $M$ to a $k$-differential form $F^*\om$ on $Q$ is defined by
		\eqn{}{(F^*\om)(q;w_1,\cd,w_k):=\int_M\rho(q,p)\om(p;V(w_1)(p),\cd,V(w_k)(p))} for $q\in Q,p\in M,w_i\in T_qQ.$
	\end{defi}
	\begin{thm}
		\eqn{}{F^*d\om=dF^*\om.}
	\end{thm}

	\B To show $F^*d\om=dF^*\om$, it is enough to apply both sides  to $(\frac{\p}{\p u_{l_1}},\cd,\frac{\p}{\p u_{l_k}})$ where $(u_1,\cd,u_k)$ is a coordinate system for $Q$. As before we set $V_i=V(\frac{\p}{\p u_i})\in\X(M).$
	\begin{eqnarray}
		& &(F^*d\om)(q;\frac{\p}{\p u_{l_1}},\cd,\frac{\p}{\p u_{l_k}})\nonumber\\&=&\int_M\rho(q)d\om(V(\frac{\p}{\p u_{l_1}}),\cd,V(\frac{\p}{\p u_{l_k}}))\nonumber\\&=&\int_M\rho(q)d\om(V_{l_1},\cd,V_{l_k})\nonumber\\&=&\sum_{i=1}^{k}(-1)^{i-1}\int_M\rho(q)V_{l_i}\om(V_{l_1},\cd,\widehat{V_{l_i}},\cd,V_{l_k})\nonumber\\&+&\sum_{i< j}(-1)^{i+j}\int_M\rho(q)\om([V_{l_i},V_{l_j}],V_{l_1},\cd,\widehat{V_{l_i}},\cd,\widehat{V_{l_j}},\cd V_{l_k})\nonumber\\&=&
		-\sum_{i=1}^{k}(-1)^{i-1}\int_M\n.(\rho V_{l_i})\om(V_{l_1},\cd,\widehat{V_{l_i}},\cd,V_{l_k})\nonumber\end{eqnarray}
	\begin{eqnarray}&+&\sum_{i< j}(-1)^{i+j}\int_M\rho(q)\om([V_{l_i},V_{l_j}],V_{l_1},\cd,\widehat{V_{l_i}},\cd,\widehat{V_{l_j}},\cd V_{l_k})\nonumber\\&=&
		\sum_{i=1}^{k}(-1)^{i-1}\int_M\frac{\p\rho}{\p u_{l_i}}\om(V_{l_1},\cd,\widehat{V_{l_i}},\cd,V_{l_k})\nonumber\\&+&\sum_{i< j}(-1)^{i+j}\int_M\rho(q)\om([V_{l_i},V_{l_j}],V_{l_1},\cd,\widehat{V_{l_i}},\cd,\widehat{V_{l_j}},\cd V_{l_k})\nonumber\\&=&
		\sum_{i=1}^{k}(-1)^{i-1}\Big(\frac{\p}{\p u_{l_i}}\int_M\rho\om(V_{l_1},\cd,\widehat{V_{l_i}},\cd,V_{l_k})\nonumber\\&-&
		\int_M\rho\frac{\p}{\p u_{l_i}}(\om(V_{l_1},\cd,\widehat{V_{l_i}},\cd,V_{l_k}))\Big)
		\nonumber\\&+&\sum_{ i< j}(-1)^{i+j}\int_M\rho(q)\om([V_{l_i},V_{l_j}],V_{l_1},\cd,\widehat{V_{l_i}},\cd,\widehat{V_{l_j}},\cd V_{l_k})\nonumber\\&=&
		\sum_{i=1}^{k}(-1)^{i-1}\Big(\frac{\p}{\p u_{l_i}}\int_M\rho\om(V_{l_1},\cd,\widehat{V_{l_i}},\cd,V_{l_k})\nonumber\\&-&
		\sum_{ j<i}(-1)^{j-1}\int_M\rho\om(\frac{\p V_{l_j}}{\p u_{l_i}},V_{l_1},\cd,\widehat{V_{l_j}},\cd,\widehat{V_{l_i}},\cd,V_{l_k})\nonumber\\&-&\sum_{ i<j}(-1)^{j-2}\int_M\rho\om(\frac{\p V_{l_j}}{\p u_{l_i}},V_{l_1},\cd,\widehat{V_{l_i}},\cd,\widehat{V_{l_j}},\cd,V_{l_k})\Big)
		\nonumber\\&+&\sum_{i< j\le }(-1)^{i+j}\int_M\rho(q)\om([V_{l_i},V_{l_j}],V_{l_1},\cd,\widehat{V_{l_i}},\cd,\widehat{V_{l_j}},\cd V_{l_k})\nonumber\\
		&=&\sum_{i=1}^{k}(-1)^{i-1}\frac{\p}{\p u_{l_i}}\int_M\rho\om(V_{l_1},\cd,\widehat{V_{l_i}},\cd,V_{l_k})\nonumber\\&+&
		\sum_{i<j}(-1)^{i+j}\int_M\rho\om(\frac{\p V_{l_j}}{\p u_{l_i}}-\frac{\p V_{l_i}}{\p u_{l_j}},V_{l_1},\cd,\widehat{V_{l_j}},\cd,\widehat{V_{l_i}},\cd,V_{l_k})
		\nonumber\\&+&\sum_{i< j\le }(-1)^{i+j}\int_M\rho(q)\om([V_{l_i},V_{l_j}],V_{l_1},\cd,\widehat{V_{l_i}},\cd,\widehat{V_{l_j}},\cd V_{l_k})\nonumber\\&=&
		\sum_{i=1}^{k}(-1)^{i-1}\frac{\p}{\p u_{l_i}}\int_M\rho\om(V_{l_1},\cd,\widehat{V_{l_i}},\cd,V_{l_k})\nonumber\\&+&
		\sum_{i<j}(-1)^{i+j}\int_M\rho\om(\frac{\p V_{l_j}}{\p u_{l_i}}-\frac{\p V_{l_i}}{\p u_{l_j}}+[V_{l_i},V_{l_j}],V_{l_1},\cd,\widehat{V_{l_j}},\cd,\widehat{V_{l_i}},\cd,V_{l_k})\nonumber\\&=&
		\sum_{i=1}^{k}(-1)^{i-1}\frac{\p}{\p u_{l_i}}\int_M\rho\om(V_{l_1},\cd,\widehat{V_{l_i}},\cd,V_{l_k})\nonumber\end{eqnarray}
	\begin{eqnarray}&=&
		\sum_{i=1}^{k}(-1)^{i-1}\frac{\p}{\p u_{l_i}}F^*\om(\frac{\p}{\p u_{l_1}},\cd,\widehat{\frac{\p}{\p u_{l_i}}},\cd,\frac{\p}{\p u_{l_k}})\nonumber\\&+&\sum_{i<j}(-1)^{i+j}F^*\om([\frac{\p}{\p u_{l_i}},\frac{\p}{\p u_{l_j}}],\frac{\p}{\p u_{l_1}},\cd,\widehat{\frac{\p}{\p u_{l_i}}},\cd,\widehat{\frac{\p}{\p u_{l_j}}},\cd \frac{\p}{\p u_{l_k}})\nonumber\\&=&dF^*\om(q;\frac{\p}{\p u_{l_1}},\cd,\frac{\p}{\p u_{l_k}}).\nonumber
	\end{eqnarray}
	\A

	Now we are going to proof a nonpure version of Stokes' theorem. 
	\begin{thm}\label{npcstokes}(\textbf{Nonpure Stokes' theorem}) For any nonpure map $F=(\rho,V)$ from a $k$-dimensional manifold to manifold $M$ and any $k$-form $\om$ over $M$
		\eqn{}{\int_QF^*d\om=\int_{\p Q}F^*\om.}
	\end{thm}
	\B 
	$\int_QF^*d\om=\int_QdF^*\om=\int_{\p Q}F^*\om.$
	\A
	
	Let us see what this theorem says in Euclidean spaces. The ordinary (pure) Stokes' theorem for a  surface $S$ in $\R^3$ parameterized
	by an ordinary (pure) map $X:D\subseteq \R^2\fl \R^3$ and for a vector field $F(x,y,z)$
	in $\R^3$ is given by  
	\eqn{}{\int_D \Big((\n\ti F)\circ X\Big).(U\ti V)dq=\int_{\partial
			D}(F\circ X).(Udu+Vdv)}
	where  $U=\frac{\p X}{\p u}, V=\frac{\p X}{\p v}$ and $dq=dudv$. Next, consider a nonpure  parameterized surface, i.e. a nonpure differentiable function $(\rho,U,V)$ form a domain $D\subseteq\R^2$ into $\R^3$, including  a scaler field $\rho=\rho(q;p)\ge0$, vanishing at infinity with respect to the variable $p$, and two vector fields $U=U(q;p)\in\R^3$,$V=V(q;p)\in\R^3$,
	$q=(u,v)\in D\subseteq\R^2, p=(x,y,z)\in\R^3$ such that   $$\frac{\partial\rho}{\partial
		u}+\nabla.(\rho U)=0,~~~\frac{\partial\rho}{\partial
		v}+\nabla.(\rho V)=0$$ and
	\eqn{con2}{\frac{\partial U}{\partial
			v}+(V.\nabla)U=\frac{\partial V}{\partial u}+(U.\nabla)V,} where $\n$ is the gradient operator with respect to the variable $p=(x,y,z)\in \R^3$. 
	Then we have
	\eqn{}{ \int_{\R^3}\int_D\rho (\n\ti F).(U\ti V)dqdp=\int_{\R^3}\int_{\partial
			D}\rho F.(Udu+Vdv)dp}
	where $dq=dudv$ and $dp=dxdydz$.

	\section{Noncommutative Analysis}
	\subsection{Noncommutative Differential Calculus}
	In this part we study  noncommutative analysis, i.e. analysis on the space of states of a noncommutative $C^*$-algebra. In the commutative case we worked with the space of nonpure states of the commutative $C^*$-algebra $C_0(M)$ 
	\eqn{}{\int_{x\in M}\rho(x)\hat{x}:C_0(M)\fl\R,~~~(\int_{x\in M}\rho(x)\hat{x})f=\int_M\rho (x)f(x)dx\in\R} where $\rho:M\fl\R$ is a  density function, i.e. it is nonnegative and $\int\rho=1$. Now we replace the commutative $C^*$-algebra with a noncommutative one $\ca $ equipped with a trace which we denote again by integration symbol 
	\eqn{}{\int:\ca \fl[0,\iy]}which plays the role of integration of the commutative case. The role of density functions of the commutative case is played with density element, i.e. self-adjoint positive elements $\rho\in \ca $ such that $\int \rho=1$.  Each density element $\rho$ induces a state on $\ca $ given by
	\eqn{}{f\in \ca \mt \int \rho f\in\R.}
	We assume that this trace is faithful, i.e. if for all $f$ we have $\int  fg=0$ we conclude that $g=0$. Moreover we  assume that the set of all elements of finite-trace  is a ideal of $\ca $ over real numbers.   
	
	\begin{defi}
		The space of  sates of $\ca $ includes all density elements, i.e. self-adjoint positive elements $\rho\in \ca $ such that $\int \rho=1$.  
	\end{defi}
	It is well-known that the integration in the commutative case satisfies in the fundamental theorem of calculus which one of its important results is the divergence theorem and thus if $\rho\in\I(M)$ and  $X\in\X(M)$ we get \eqn{}{\int\n.(\rho X)=0.} The replacement for this formula in the noncommutative case is the following  property of trace
	\eqn{}{\int[\rho,X]=0,~~~~~~\fa \rho,X\in\ca.}Thus it seems that the concept of trace in noncommutative case is a suitable replacement for the concept of integration in commutative case. 
	
	After introducing the space of nonpure states which is the replacement for the space of nonpure states of a commutative	algebra $C_0(M)$, which includes all density functions we are ready for the next important step which is to establish differential calculus on this space. Recall that in the commutative case $C_0(M)$ following the pure case we said a path $\rho_t:M\fl[0,\iy)$ of nonpure sates is differentiable if there exists a time-dependent vector field $V_t$ over $M$ such that
	\eqn{1}{\frac{d}{dt}\int\rho f=\int\rho Vf} for all smooth functions $f$ on $M$. Here $Vf$ is the derivative of $f$ in direction of $V$. Now since each vector field $V$ gives a derivation on the algebra $C^\iy(M)$ we must replace the vector field $V_t$ with a derivation on the algebra $\ca $ in the noncommutative case, i.e. with the mapping
	\eqn{}{f\in \ca \mt[V,f].}
	\begin{defi}
		A path $\rho_t$ of states, i.e. density elements, is called differentiable if there exists a path $V_t\in \ca $ such that
		\eqn{fracddtintM}{\frac{d}{dt}\int \rho_t f =\int \rho_t[V_t,f ] } for all $f \in \ca $.
	\end{defi}
	In the commutative case by the aid of divergence theorem we could prove that the equation (\ref{1}) is equivalent with the following equation
	\eqn{}{\frac{\p\rho}{\p t}+\n.(\rho V)=0.}
	Now in the noncommutative case we have the following fact.
	\begin{thm}
		The equation (\ref{fracddtintM}) is equivalent with 
		\eqn{}{\frac{\p \rho }{\p t}+[V,\rho ]=0.}
	\end{thm}
	\B Since the trace is linear we have
	$\frac{d}{dt}\int \rho_t f =\int\frac{\p \rho }{\p t}f $. On the other hand $\int \rho [V,f ]=\int[V,\rho f ]-\int[V,\rho ]f =-\int[V,\rho ]f .$ Thus $\int(\frac{\p \rho }{\p t}+[V,\rho ])f =0$ for all $f $. Thus by the faithfulness of the trace we conclude the desired result. The converse is similarly proved. \A
	
	In the nonpure commutative case we called a two-parameter density function $\rho(t,s)$  differentiable if there exist two $(t,s)$-dependent vector fields $V$ and $W$ over $M$ satisfying
	$\frac{\p\rho}{\p t}+\n.(\rho V)=0$ and $\frac{\p\rho}{\p s}+\n.(\rho W)=0$.
	\begin{defi}
		A two-parameter density element $\rho (t,s)\in \ca $ is called differentiable if there exist two $(t,s)$-dependent elements $V$ and $W$ in $\ca $ satisfying
		\eqn{fracpm}{\frac{\p \rho }{\p t}+[V,\rho ]=0,~~~~~\frac{\p \rho }{\p s}+[W,\rho ]=0.}
	\end{defi}
	In the pure commutative case $C(M)$ we can prove the identity
	\eqn{}{\frac{\p v}{\p s}-\frac{\p w}{\p t}=0,} where $v(t,s)=\frac{\p x}{\p t}$ and $w(t,s)=\frac{\p x}{\p s}$ for some smooth pure function $x(t,s)\in M$. The analogue of this identity for the nonpure commutative case is 
	\eqn{nbigrho}{\n.\Big(\rho(\frac{\p V}{\p s}-\frac{\p W}{\p t}-[V,W])\Big)=0.}
	Here $[V,W]$ is the Lie bracket of two vector fields. 
	The analogue of this identity for the nonpure noncommutative case is the following fact.
	\begin{thm}
		\eqn{}{[\frac{\p V}{\p s}-\frac{\p W}{\p t}-[V,W],\rho ]=0.}Here $[,]$ is the commutator of two elements.
	\end{thm}
	\B First proof: we have
	$\frac{\p^2\rho }{\p s\p t}=\frac{\p}{\p s}[\rho ,V]=[\frac{\p \rho }{\p s},V]+[\rho ,\frac{\p V}{\p s}]=[[\rho ,W],V]+[\rho ,\frac{\p V}{\p s}].$ Thus using the Jacobi identity in the forth line below we have
	\begin{eqnarray}0&=&\frac{\p^2\rho }{\p s\p t}-\frac{\p^2\rho }{\p t\p s}\nonumber\\&=&[[\rho ,W],V]-[[\rho ,V],W]+[\rho ,\frac{\p V}{\p s}]-[\rho ,\frac{\p W}{\p t}]\nonumber\\
		&=&[[\rho ,W],V]+[[V,\rho ],W]+[\rho ,\frac{\p V}{\p s}]-[\rho ,\frac{\p W}{\p t}]\nonumber\\
		&=&-[[W,V],\rho ]+[\rho ,\frac{\p V}{\p s}]-[\rho ,\frac{\p W}{\p t}]\nonumber\\
		&=&[\rho ,\frac{\p V}{\p s}-\frac{\p W}{\p t}-[V,W]].\nonumber\end{eqnarray}Second proof: by the formula  (\ref{fracddtintM}) we have for all $f $
	\begin{eqnarray}\frac{\p^2}{\p s\p t}\int \rho f &=&
		\int(\frac{\p \rho }{\p s}[V,f ]+\rho [\frac{\p V}{\p s},f ])\nonumber\\
		&=&\int([\rho ,W][V,f ]+\rho [\frac{\p V}{\p s},f ])\nonumber\\
		&=&\int[V,[\rho ,W]f ]-\int[V,[\rho ,W]]f +\int[\frac{\p V}{\p s},\rho f ]-\int[\frac{\p V}{\p s},\rho ]f \nonumber\\
		&=&\int([[\rho ,W],V]+[\rho ,\frac{\p V}{\p s}])f .\nonumber
	\end{eqnarray}
	Hence $0=\frac{\p^2}{\p s\p t}\int \rho f -\frac{\p^2}{\p t\p s}\int \rho f =\int([[\rho ,W],V]+[\rho ,\frac{\p V}{\p s}]-[[\rho ,W],V]-[\rho ,\frac{\p V}{\p s}])f $. Thus $[[\rho ,W],V]+[\rho ,\frac{\p V}{\p s}]-[[\rho ,V],W]-[\rho ,\frac{\p W}{\p t}]=0.$ Now the remainder of the proof is similar to the previous proof.
	\A
	
	In the nonpure commutative case we assumed the stronger version of (\ref{nbigrho}), i.e. we assumed that
	\eqn{}{\frac{\p V}{\p s}-\frac{\p W}{\p t}=[V,W].}
	Similarly we add the following hypothesis to the definition of differentiability in the noncommutative case.
	\begin{defi}
		A two-parameter density element $\rho (t,s)\in \ca $ is called differentiable if there exist two $(t,s)$-dependent elements $V$ and $W$ in $\ca $ satisfying (\ref{fracpm}) and 
		\eqn{}{\frac{\p V}{\p s}-\frac{\p W}{\p t}=[V,W].}
	\end{defi}
	More generally we can define differentiability of maps from a manifold to the space $\s(\ca )$.
	\begin{defi}
		A \textbf{differentiable    function} from open subset
		$U\se\R^m$ into the space $\s(\ca )$ of states of the algebra  $A$ is a system
		$(\rho ,V_1,\cd,V_m)$ including a   function
		\eqn{}{\rho:U\fl \s(\ca )}which is differentiable too, together
		with differentiable maps \eqn{}{V_j:U\fl \ca } $1\le
		j\le m,$ such that \eqn{}{\frac{\p \rho }{\p u_j}+[V_j,\rho ]=0}
		and \eqn{}{\frac{\p V_i}{\p u_j}-\frac{\p V_j}{\p u_i}=[V_i,V_j].}
	\end{defi}
	
	\begin{defi}\label{smoouncerfun}
		A \textbf{differentiable    function} from manifold $P$ into the space $\s(\ca )$ of states of the algebra  $\ca $ is a system $(\rho,V)$
		including a   map \eqn{}{\rho :P\fl \s(\ca )}which is differentiable, together
		with a differentiable linear map on each fiber 
		\eqn{}{V:TP\fl \ca } 
		such that  for any coordinate system
		$u=(u_1,\cd,u_m)\in U\se\R^m$ for $P$ \eqn{pmpuj}{\frac{\p \rho }{\p
				u_j}+[V_j,\rho ]=0} where $V_j:U\fl
		\ca ,V_j(u):=V(\frac{\p}{\p u_j}(u))$ and
		\eqn{fracpvipuj}{\frac{\p V_i}{\p u_j}-\frac{\p
				V_j}{\p u_i}=[V_i,V_j].}
		
	\end{defi}
	\begin{prop}\label{dprho}
		The  equations (\ref{pmpuj}) and (\ref{fracpvipuj}) do not depend to coordinate systems.
	\end{prop}
	\B Let $B_i(v):=V_{p}(\frac{\p}{\p v_i}(p))$ where $v$ is the
	coordinates of point $p$  then since $\frac{\p}{\p
		v_i}=\sum_j\frac{\p u_j}{\p v_i}\frac{\p}{\p u_j}$ we get
	$B_i(v)=\sum_j\frac{\p u_j}{\p v_i}V_j(u)$. Thus
	\begin{eqnarray}& &\frac{\p \rho (v)}{\p v_i}+[B_i(v),\rho (v)]
		\nonumber\\&=& \sum_j\frac{\p  u_j}{\p v_i}\frac{\p \rho (u)}{\p
			u_j}+\sum_j\frac{\p u_j}{\p
			v_i}[V_j(u),\rho (u)]\nonumber\\&=&\sum_j\frac{\p u_j}{\p
			v_i}[\frac{\p \rho (u)}{\p u_j},[V_j(u),\rho (u)].\nonumber
	\end{eqnarray}Thus since the matrix $(\frac{\p u_j}{\p
		v_i})$ is invertible we conclude that $\frac{\p \rho (v)}{\p
		v_i}+[B_i(v),\rho (v)]=0$ if and only if
	$\frac{\p \rho (u)}{\p u_i}+[V_i(u),\rho (u)]=0.$
	
	Next,  
	$$\frac{\p B_i}{\p v_j}=\sum_k\frac{\p^2u_k}{\p v_j\p
		v_i}V_k+\sum_{kl}\frac{\p u_k}{\p v_i}\frac{\p u_l}{\p
		v_j}\frac{\p V_k}{\p u_l}$$ Thus $$\frac{\p B_i}{\p v_j}-\frac{\p
		B_j}{\p v_i}=\sum_{kl}\frac{\p u_k}{\p v_i}\frac{\p u_l}{\p
		v_j}(\frac{\p V_k}{\p u_l}-\frac{\p V_l}{\p u_k})$$ and
	$$[B_i,B_j]=\sum_{kl}\frac{\p u_k}{\p v_i}\frac{\p
		u_l}{\p v_j}[V_k,V_l].$$ Hence
	$$\frac{\p B_i}{\p u_j}-\frac{\p B_j}{\p u_i}-[B_i,B_j]=
	\sum_{kl}\frac{\p u_k}{\p v_i}\frac{\p u_l}{\p v_j}(\frac{\p
		V_k}{\p u_l}-\frac{\p V_l}{\p u_k}-[V_k,V_l]).$$ Thus since the
	matrix $(\p
	u_j/\p v_i)$ is invertible we conclude that\\
	$$\frac{\p B_i}{\p u_j}-\frac{\p B_j}{\p
		u_i}-[B_i,B_j]=0$$ if and only if $$\frac{\p V_k}{\p u_l}-\frac{\p
		V_l}{\p u_k}-[V_k,V_l]=0.$$ \A
	
	\begin{rem}\label{2}
		By comparing the above theory in the commutative and noncommutative case we realize that the theory in  commutative case is not just a special case of the theory in noncommutative case. In fact if we apply the theory in noncommutative case to the commutative algebras we get trivial results since for example a differentiable nonpure curve becomes just a constant curve. So the theories in commutative case and noncommutative cases are independent of each other but parallel with each other. Namely there is a deep analogy between them as follows. In fact there is a dictionary between them as follows which translate any concept, proposition and proof from one side to the other side:\\ 
		\[density ~functions~\Longleftrightarrow~ density~ operators,\]
		\[vector~ fields~\Longleftrightarrow~ arbitrary~ ~operators~ in ~\ca,\] 
		\[directional~ derivative ~Vf~\Longleftrightarrow~  ~commutator~[V,f]~,\]
		\[Lie ~bracket~[V,W]~\Longleftrightarrow~commutator ~[V,W],\] 
		\[integration ~\Longleftrightarrow~  trace,\] 
		\[divergence~ \n.(fV)~\Longleftrightarrow~commutator ~[V,f],\] 
	\end{rem}
	\subsection{Noncommutative Integral Calculus}
	In this section, we state and prove the  noncommutative version of the Stokes theorem. We first introduce the noncommutative version of the differential forms and de Rham cohomology.
	\begin{defi}
		The $k$-differential forms over algebra $\ca$ is defined to be the linear maps $\om:\Lan^k\ca\fl\ca$ where $\Lan^k\ca$ is the $k$-th exterior algebra. In other words a $k$-differential form on $\ca$ is a $k$-linear alternating map $\om:\ca^{\otimes k}\fl\ca$. The space of  $k$-differential forms on $\ca$ is denoted by $\Om^k(\ca)$. The differential operator $d:\Om^k(\ca)\fl\Om^{k+1}(\ca)$ is the  Chevalley-Eilenberg operator. The de Rham cohomology of the algebra $\ca$ is defined to be the Chevalley-Eilenberg cohomology of the Lie algebra $(A,[,])$ with coefficients in the module $A$ with adjoint action i.e. the bracket, of $A$ on itself. 
	\end{defi} 
	
	\begin{defi}
		Let $F=(\rho,V)$ be a nonpure map from oriented manifold $Q$ to algebra $\ca$, i.e. $\rho:Q\fl \s(\ca), V:TQ\fl \ca$. The pullback of $k$-differential form $\om$ on $\ca$ to a $k$-differential form $F^*\om$ on $Q$ is defined by
		\eqn{}{(F^*\om)(q;w_1,\cd,w_k):=\int\rho(q)\om(V(w_1)\otimes\cd\otimes V(w_k))} for $q\in Q,w_i\in T_qQ.$
	\end{defi}
	\begin{thm}\label{npcstokes}(\textbf{Nonpure commutative Stokes theorem}) 
		We have $F^*d\om=dF^*\om$ and 
		\eqn{}{\int_QF^*d\om=\int_{\p Q}F^*\om.}
	\end{thm}
	\B We first prove that $F^*d\om=dF^*\om$. it is enough to apply both sides  to $(\frac{\p}{\p u_{l_1}},\cd,\frac{\p}{\p u_{l_k}})$ where $(u_1,\cd,u_k)$ is a coordinate system for $Q$. As before we set $V_i=V(\frac{\p}{\p u_i})\in\ca.$
	\begin{eqnarray}
		& &(F^*d\om)(q;\frac{\p}{\p u_{l_1}},\cd,\frac{\p}{\p u_{l_k}})\nonumber\\&=&\int \rho(q)d\om(V(\frac{\p}{\p u_{l_1}}),\cd,V(\frac{\p}{\p u_{l_k}}))\nonumber\\&=&\int \rho(q)d\om(V_{l_1},\cd,V_{l_k})\nonumber\\&=&\sum_{i=1}^{k}(-1)^{i-1}\int \rho(q)[V_{l_i},\om(V_{l_1},\cd,\widehat{V_{l_i}},\cd,V_{l_k})]\nonumber\\&+&\sum_{i< j}(-1)^{i+j}\int \rho(q)\om([V_{l_i},V_{l_j}],V_{l_1},\cd,\widehat{V_{l_i}},\cd,\widehat{V_{l_j}},\cd V_{l_k})\nonumber\\&=&
		-\sum_{i=1}^{k}(-1)^{i-1}\int [V_{l_i},\rho ]\om(V_{l_1},\cd,\widehat{V_{l_i}},\cd,V_{l_k})\nonumber\\&+&\sum_{i< j}(-1)^{i+j}\int \rho(q)\om([V_{l_i},V_{l_j}],V_{l_1},\cd,\widehat{V_{l_i}},\cd,\widehat{V_{l_j}},\cd V_{l_k})\nonumber\\&=&
		\sum_{i=1}^{k}(-1)^{i-1}\int \frac{\p\rho}{\p u_{l_i}}\om(V_{l_1},\cd,\widehat{V_{l_i}},\cd,V_{l_k})\nonumber\\&+&\sum_{i< j}(-1)^{i+j}\int \rho(q)\om([V_{l_i},V_{l_j}],V_{l_1},\cd,\widehat{V_{l_i}},\cd,\widehat{V_{l_j}},\cd V_{l_k})\nonumber\\&=&
		\sum_{i=1}^{k}(-1)^{i-1}\Big(\frac{\p}{\p u_{l_i}}\int \rho\om(V_{l_1},\cd,\widehat{V_{l_i}},\cd,V_{l_k})\nonumber\\&-&
		\int \rho\frac{\p}{\p u_{l_i}}(\om(V_{l_1},\cd,\widehat{V_{l_i}},\cd,V_{l_k}))\Big)
		\nonumber\\&+&\sum_{ i< j}(-1)^{i+j}\int \rho(q)\om([V_{l_i},V_{l_j}],V_{l_1},\cd,\widehat{V_{l_i}},\cd,\widehat{V_{l_j}},\cd V_{l_k})\nonumber\\&=&
		\sum_{i=1}^{k}(-1)^{i-1}\Big(\frac{\p}{\p u_{l_i}}\int \rho\om(V_{l_1},\cd,\widehat{V_{l_i}},\cd,V_{l_k})\nonumber\\&-&
		\sum_{ j<i}(-1)^{j-1}\int \rho\om(\frac{\p V_{l_j}}{\p u_{l_i}},V_{l_1},\cd,\widehat{V_{l_j}},\cd,\widehat{V_{l_i}},\cd,V_{l_k})\nonumber\\&-&\sum_{ i<j}(-1)^{j-2}\int \rho\om(\frac{\p V_{l_j}}{\p u_{l_i}},V_{l_1},\cd,\widehat{V_{l_i}},\cd,\widehat{V_{l_j}},\cd,V_{l_k})\Big)
		\nonumber\\&+&\sum_{i< j\le }(-1)^{i+j}\int \rho(q)\om([V_{l_i},V_{l_j}],V_{l_1},\cd,\widehat{V_{l_i}},\cd,\widehat{V_{l_j}},\cd V_{l_k})\nonumber
	\end{eqnarray}
	\begin{eqnarray}
		&=&\sum_{i=1}^{k}(-1)^{i-1}\frac{\p}{\p u_{l_i}}\int \rho\om(V_{l_1},\cd,\widehat{V_{l_i}},\cd,V_{l_k})\nonumber\\&+&
		\sum_{i<j}(-1)^{i+j}\int \rho\om(\frac{\p V_{l_j}}{\p u_{l_i}}-\frac{\p V_{l_i}}{\p u_{l_j}},V_{l_1},\cd,\widehat{V_{l_j}},\cd,\widehat{V_{l_i}},\cd,V_{l_k})
		\nonumber\\&+&\sum_{i< j\le }(-1)^{i+j}\int \rho(q)\om([V_{l_i},V_{l_j}],V_{l_1},\cd,\widehat{V_{l_i}},\cd,\widehat{V_{l_j}},\cd V_{l_k})\nonumber\\&=&
		\sum_{i=1}^{k}(-1)^{i-1}\frac{\p}{\p u_{l_i}}\int \rho\om(V_{l_1},\cd,\widehat{V_{l_i}},\cd,V_{l_k})\nonumber\\&+&
		\sum_{i<j}(-1)^{i+j}\int \rho\om(\frac{\p V_{l_j}}{\p u_{l_i}}-\frac{\p V_{l_i}}{\p u_{l_j}}+[V_{l_i},V_{l_j}],V_{l_1},\cd,\widehat{V_{l_j}},\cd,\widehat{V_{l_i}},\cd,V_{l_k})\nonumber\\&=&
		\sum_{i=1}^{k}(-1)^{i-1}\frac{\p}{\p u_{l_i}}\int \rho\om(V_{l_1},\cd,\widehat{V_{l_i}},\cd,V_{l_k})\nonumber\\&=&
		\sum_{i=1}^{k}(-1)^{i-1}\frac{\p}{\p u_{l_i}}F^*\om(\frac{\p}{\p u_{l_1}},\cd,\widehat{\frac{\p}{\p u_{l_i}}},\cd,\frac{\p}{\p u_{l_k}})\nonumber\\&+&\sum_{i<j}(-1)^{i+j}F^*\om([\frac{\p}{\p u_{l_i}},\frac{\p}{\p u_{l_j}}],\frac{\p}{\p u_{l_1}},\cd,\widehat{\frac{\p}{\p u_{l_i}}},\cd,\widehat{\frac{\p}{\p u_{l_j}}},\cd \frac{\p}{\p u_{l_k}})\nonumber\\&=&dF^*\om(q;\frac{\p}{\p u_{l_1}},\cd,\frac{\p}{\p u_{l_k}}).\nonumber
	\end{eqnarray}Next we have 
	$\int_QF^*d\om=\int_QdF^*\om=\int_{\p Q}F^*\om.$ In fact the above proof is just the translation of the proof in nonpure commutative case under the dictionary given in Remark \ref{2}. 
	\A


\begin{thebibliography}{widestlabel}
		\bibitem{AG} Ambrosio, L., Gigli, N.: Construction of the parallel transport in the Wasserstein space. Methods Appl. Anal. (1) 15,  1–30 (2008)
		\bibitem{AGS} 
		Ambrosio, L., Gigli, N., Savaré, G.: Gradient Flows
		in Metric Spaces and in the Space
		of Probability Measures. Birkhäuser Verlag, (2005)
		\bibitem{Con}
		Connes, C.: Noncommutative Geometry. Academic Press, (1994)
		\bibitem{F}
		Franco, N.: The Lorentzian distance formula in noncommutative geometry. In Non-Regular Spacetime Geometry. J. Phys. Conf. Ser. 968 (2018) 
		\bibitem{GNT}
		Gangbo, W., Nguyen, T., Tudorascu, A.: 	Hamilton-Jacobi Equations in the Wasserstein Space. Methods Appl. Anal. 
		(2) 15, 155–184 (2008)
		\bibitem{L1} 
		Lott, J.: Some Geometric Calculations on Wasserstein Space. Commun. Math. Phys. (2) 277, 423–437  (2008)
		\bibitem{L2} 
		Lott, J.: An Intrinsic Parallel Transport in Wasserstein Space. Proc. Am. Math. Soc. (12) 145, 5329–5340  (2017)
		
		\bibitem{LV} 
		Lott, J., Villani, C.: Ricci curvature for metric-measure spaces via optimal transport. Ann. Math. (3) 169, 903-991  (2009)
		\bibitem{M1} 
		Martinetti, P.: From Monge to Higgs: a Survey of Distance Computations in Noncommutative
		Geometry. In Workshop on
		Noncommutative Geometry and Optimal Transport (Besançon, France, November 27, 2014), Contemp. Math.  676 (2016)
		\bibitem{O} 
		Otto, F.: The geometry of dissipative evolution equations: the porous meduim equation. Comm.
		Partial Differential Equations. 26, 101–174  (2001)
		\bibitem{R1}  Rieffel, M. A.: Metrics on states from actions of compact groups. Doc. Math. 3,
		215–229  (1998) 
		\bibitem{R2}  Rieffel, M. A.: Metrics on state spaces. Doc. Math. 4, 559–600  (1999).
		\bibitem{V} Villani, C.: Optimal Transport,
		Old and New. Springer (2009)
	\end{thebibliography}
\end{document}